\documentclass[a4paper,12pt]{article}
\usepackage{amsthm,amsmath}
\usepackage{amssymb}
\usepackage[latin1]{inputenc}
\usepackage[all]{xypic}
\usepackage{enumerate}
\usepackage{a4wide}

\newtheorem{theorem}{Theorem}[section]
\newtheorem{proposition}[theorem]{Proposition}

\newtheorem{lemma}[theorem]{Lemma}
\newtheorem*{theorem*}{Theorem}
\newtheorem*{proposition*}{Proposition}
\newtheorem*{corollary*}{Corollary}
\newtheorem*{lemma*}{Lemma}
\theoremstyle{definition}

\newtheorem{remark}[theorem]{Remark}
\newtheorem*{remark*}{Remark}

\newtheorem*{definition*}{Definition}

\newcommand{\coring}[1]{\mathfrak{#1}}
\newcommand{\tensor}[1]{\otimes_{#1}}
\newcommand{\tensfun}[1]{\underset{{#1}}{\otimes}}

\newcommand{\rcomod}[1]{\mathcal{M}^{#1}}
\newcommand{\rmod}[1]{\mathcal{M}_{#1}}
\newcommand{\lmod}[1]{{}_{#1}\mathcal{M}}

\newcommand{\cotensor}[1]{\square_{#1}}
\newcommand{\lcomod}[1]{{}^{#1}\mathcal{M}}

\renewcommand{\hom}[3]{\mathrm{Hom}_{#1}(#2,#3)}
\newcommand{\End}[2]{\mathrm{End}_{#1-cor}(#2)}
\newcommand{\Aut}[2]{\mathrm{Aut}_{#1-cor}(#2)}
\newcommand{\rend}[2]{\mathrm{End}({#2}_{#1})}
\newcommand{\lend}[2]{\mathrm{End}({}_{#1}#2)}

\newcommand{\rcomatrix}[2]{#2^* \tensor{#1} #2}

\newcommand{\swe}[3]{#1 \tensor{#2} #3}
\newcommand{\ibs}{\mathbf{I}_B(S)}
\newcommand{\ilbs}{\mathbf{I}_B^l(S)}
\newcommand{\irbs}{\mathbf{I}_B^r(S)}
\newcommand{\invbs}{\mathbf{Inv}_B(S)}

\begin{document}
\title{Comatrix Corings and Invertible Bimodules\footnote{Research supported by the grant BFM2001-3141 from the Ministerio de Ciencia y Tecnolog\'{\i}a of Spain}}
\author{L. El Kaoutit\footnote{Supported by
the grant SB2003-0064 from the Mi\-nis\-terio de Educaci{\'o}n,
Cultura y Deporte of Spain. } { \normalsize and } J.
G{\'o}mez-Torrecillas \\ \normalsize \footnotesize{Departamento de
{\'A}lgebra, Universidad de Granada, E18071 Granada, Spain.}
\\ \normalsize e-mails: \textsf{kaoutit@ugr.es}\,\, \textsf{gomezj@ugr.es} }

\date{December 31, 2004}

\maketitle

\begin{abstract}
We extend Masuoka's Theorem \cite{Masuoka:1989} concerning the
isomorphism between the group of invertible bimodules in a
non-commutative ring extension and the group of automorphisms of
the associated Sweedler's canonical coring, to the class of finite
comatrix corings introduced in \cite{Kaoutit/Gomez:2003a}.
\end{abstract}

\section*{Introduction}

Comatrix corings were introduced by the authors in
\cite{Kaoutit/Gomez:2003a} to give a structure theorem of all
cosemisimple corings. This construction generalizes Sweedler's
canonical corings \cite{Sweedler:1975}, and provides a version of
descent theory for modules \cite[Theorem
3.10]{Kaoutit/Gomez:2003a}. Sweedler's canonical corings and their
automorphisms were the key tool in \cite{Masuoka:1989} to give a
non-commutative version of the fact that the relative Picard group
attached to any commutative ring extension is isomorphic to the
Amistur 1-cohomology for the units-functor due to Grothendieck's
faithfully flat descent.

In this note we extend, by using different methods, the main
result of \cite[$\S$2]{Masuoka:1989} to the context of comatrix
corings. In fact, we apply ideas and recent results from
\cite{Gomez:2002} and \cite{Kaoutit/Gomez:2003a}, and the present
paper can be already seen as natural continuation of the theory
developed in \cite{Kaoutit/Gomez:2003a}.

The first section is rather technical, and it is devoted to prove
that there is an adjoint pair of functors between the category of
comodules over a given comatrix coring and the category of
comodules over its associated Sweedler's canonical coring. This
adjunction will have a role in the proof of the main result.
Section 2 is the core of the paper, as it contains the
aforementioned isomorphism of groups (Theorem \ref{resultado-1}).
The maps connecting bimodules and coring automorphisms are at a
first glance different than the maps constructed in
\cite{Masuoka:1989}. However, they are neatly related, as
Proposition \ref{(b)1} shows.

All rings considered in this note are algebras with $1$ over
commutative ground base ring $K$. A right or left module, means a
unital module. All bimodules over rings are central
$K$--bimodules. If $A$ is any ring, then we denote by $\rmod{A}$
(resp. $\lmod{A}$) the category of all right (resp. left)
$A$--modules. The opposite ring of $A$ will be denoted by $A^o$,
its multiplication is defined by $a_2^o a_1^o= (a_1a_2)^o$,
$a_1^o, a_2^o \in A^o$ (i.e. $a_1, a_2 \in A$). As usual, some
special convention will be understood for the case of
endomorphisms rings of modules. Thus, if $X_A$ is an object of
$\rmod{A}$, then its endomorphisms ring will be denoted by
$\rend{A}{X}$, while if ${}_AY$ is left $A$--module, then its
endomorphisms ring, denoted by $\lend{A}{Y}$, is, by definition,
the opposite of the endomorphisms ring of $Y$ as an object of the
category $\lmod{A}$. In this way $X$ is an
$(\rend{A}{X},A)$--bimodule, while $Y$ is an
$(A,\lend{A}{Y})$--bimodule. The opposite left $A^o$--module of
$X_A$, will be denoted by $X^o$, the action is given by $a^o x^o=
(xa)^o$, $a^o \in A^o$, $x^o \in X^o$. Of course, if $f : X
\rightarrow W$ is right $A$--linear map, then its opposite map
$f^o: X^o \rightarrow W^o$ is left $A^o$--linear which is defined
by $f^o(x^o)=(f(x))^o$, for all $x^o \in X^o$. The same process
will be applied on bimodules and bilinear maps. For any
$(B,A)$--bimodule $M$ we denote by $M^* = \mathrm{Hom}(M_A,A_A)$
its right dual and by ${}^*M  =  \mathrm{Hom}({}_BM,{}_BB)$ its
left dual. $M^*$ and ${}^*M$ are considered, in a natural way, as
an $(A,B)$--bimodules.

Recall from \cite{Sweedler:1975} that an $A$--coring is a
three-tuple
$(\coring{C},\Delta_{\coring{C}},\varepsilon_{\coring{C}})$
consisting of an $A$--bimodule $\coring{C}$ and the two
$A$--bilinear maps
$$\xymatrix@C=50pt{\coring{C} \ar@{->}^-{\Delta_{\coring{C}}}[r] & \coring{C}\tensor{A}\coring{C}},\quad
\xymatrix@C=30pt{ \coring{C}
\ar@{->}^-{\varepsilon_{\coring{C}}}[r] & A}$$ such that
$(\Delta_{\coring{C}}\tensor{A}\coring{C}) \circ
\Delta_{\coring{C}} = (\coring{C}\tensor{A}\Delta_{\coring{C}})
\circ \Delta_{\coring{C}}$ and
$(\varepsilon_{\coring{C}}\tensor{A}\coring{C}) \circ
\Delta_{\coring{C}}=(\coring{C}\tensor{A}\varepsilon_{\coring{C}})
\circ \Delta_{\coring{C}}= \coring{C}$. A morphism of an
$A$--corings, is an $A$--bilinear map $\phi: \coring{C}
\rightarrow \coring{D}$ which satisfies: $\varepsilon_{\coring{D}}
\circ \phi = \varepsilon_{\coring{C}}$ and $\Delta_{\coring{D}}
\circ \phi = (\phi \tensor{A} \phi) \circ \Delta_{\coring{C}}$. A
right $\coring{C}$--comodule is a pair $(M,\rho_{M})$ consisting
of right $A$--module and a right $A$--linear map $\rho_{M}: M
\rightarrow M\tensor{A}\coring{C}$, called right
$\coring{C}$--coaction, such that
$(M\tensor{A}\Delta_{\coring{C}}) \circ \rho_M =
(\rho_M\tensor{A}\coring{C}) \circ \rho_M$ and
$(M\tensor{A}\varepsilon_{\coring{C}}) \circ \rho_M=M$. Left
$\coring{C}$--comodules are symmetrically defined, and we will use
the Greek letter $\lambda_{-}$ to denote theirs coactions. For
more details on comodules, definitions and basic properties of
bicomodules and the cotensor product, the reader is referred to
\cite{Brzezinski/Wisbauer:2003} and its bibliograpy.

\section{Comatrix coring and adjunctions}\label{Sect1}

Throughout this section $\Sigma$ will be a fixed $(B,A)$--bimodule
which is finitely generated and projective as right $A$--module
with a fixed dual basis $\{(e_i,e_i^*)\}_{1 \leq i \leq n} \subset
\Sigma \times \Sigma^*$.  Let $S = \rend{A}{\Sigma}$ its right
endomorphisms ring, and let $\lambda : B \rightarrow S$ be the
canonical associated ring extension. It is known that there is a
$S$--bimodule isomorphism
\begin{equation}\label{xi}
\xymatrix@R=0pt{\xi: \Sigma \tensor{A} \Sigma^* \ar@{->}[r] & S =
\rend{A}{\Sigma}
\\ u \tensor{A} v^* \ar@{|->}[r] & [x \mapsto u v^*(x)] \\
\sum_i e_i \tensor{A} e_i^*s=\sum_i se_i \tensor{A} e_i^* & s
\ar@{|->}[l] }
\end{equation}
With this identification the product of $S$ (the composition)
satisfies
\begin{equation}\label{tensor}
\begin{array}{c}
s(u \tensor{A} u^*) = s(u) \tensor{A} u^*, \\ (u \tensor{A} u^*) s
= u \tensor{A} u^* s, \\ (u \tensor{A} u^*) (v \tensor{A} v^*) = u
u^*(v) \tensor{A} v^* = u \tensor{A} u^*(v) v^*,
\end{array}
\end{equation}
for every $s \in S$, $u, v \in \Sigma$, $v^*, u^* \in \Sigma^*$.
By \cite[Proposition 2.1]{Kaoutit/Gomez:2003a}, the $A$--bimodule
$\rcomatrix{B}{\Sigma}$ is an $A$--coring with the following
comultiplication and counit
\[
\Delta_{\rcomatrix{B}{\Sigma}}(u^* \tensor{B}u) = \sum_i u^*
\tensor{B} e_i \tensor{A} e_i^* \tensor{B} u,\quad
\varepsilon_{\rcomatrix{B}{\Sigma}}(u^* \tensor{B}u)= u^*(u).
\]
The map $\Delta_{\rcomatrix{B}{\Sigma}}$ is independent on the
choice of the right dual basis of $\Sigma_A$, see \cite[Remark
2.2]{Kaoutit/Gomez:2003a}. This coring is known as \emph{ the
comatrix coring} associated to the bimodule $\Sigma$.

\begin{remark}\label{lcomcor}
One can define a comatrix coring using a bimodule which is a
finitely generated and projective left module. However, the
resulting coring is isomorphic to the comatrix coring defined by
the left dual module. To see this, consider ${}_A\Lambda_B$ any
bimodule such that ${}_A\Lambda$ is finitely generated and
projective module with a fixed left dual basis
$\{f_j,{}^*f_j\}_j$. Put ${}_B\Sigma_A ={}_B{}^*\Lambda_A$, the
set $\{{}^*f_j,f_j^*\}_j$ where $f_j^* \in \Sigma^*$ are defined
by $f_j^*(u)=u(f_j)$, for all $u \in \Sigma$ and $j$; form a right
dual basis for $\Sigma_A$. The isomorphism of corings is given by
$$ \xymatrix@R=0pt{ \rcomatrix{B}{\Sigma} \ar@{->}^-{\cong}[r] &
\Lambda\tensor{B}{}^*\Lambda \\ u^* \tensor{B}{}^*v \ar@{|->}[r] &
(\sum_ju^*({}^*f_j)f_j)\tensor{B}{}^*v }$$ The proof is direct,
using the above duals basis, and we leave it to the reader.
\end{remark}

Keeping the notations before the Remark \ref{lcomcor}, we have
that the right (resp. left) $A$--module $\Sigma$ (resp.
$\Sigma^*$) is a right (resp. left)
$\rcomatrix{B}{\Sigma}$--comodule with left (resp. right)
$B$--linear coaction: $$\rho_{\Sigma} : \Sigma \longrightarrow
\Sigma \tensor{A}\rcomatrix{A}{\Sigma},\quad (u \mapsto \sum_i e_i
\tensor{A} e_i^* \tensor{B} u),$$ for every $u \in \Sigma$, and
$$\lambda_{\Sigma^*} : \Sigma^* \longrightarrow \rcomatrix{A}{\Sigma}
\tensor{A} \Sigma^*,\quad (u^* \mapsto \sum_i u^* \tensor{A} e_i
\tensor{B} e_i^*),$$ for every $u^* \in \Sigma^*$. Furthermore,
the natural right $A$--linear isomorphism $\Sigma \cong
{}^*(\Sigma^*)$ turns out to be a right
$\rcomatrix{B}{\Sigma}$--colinear isomorphism. Associated to the
ring extension $\lambda: B \rightarrow S$, we consider also the
canonical Sweedler $S$--coring $\swe{S}{B}{S}$ whose
comultiplication is given by
$\Delta_{S\tensor{B}S}(s\tensor{B}s')=s\tensor{B}1\tensor{S}1\tensor{B}s'$,
$s,s' \in S$, and the counit is the usual multiplication.

\bigskip
The aim of this section is to establish an adjunction between the
category of right $\rcomatrix{B}{\Sigma}$--comodules and the
category of right $\swe{S}{B}{S}$--comodules. Recall first that
this last category is isomorphic to the category of descent data
associated to the extension $B \rightarrow S$, (cf.
\cite{Nuss:1997}, \cite{Brzezinski:2002}). This isomorphism of
categories will be implicitly used in the sequel. For every left
$S$--module $Y$ and right $S$--module $Z$, we denote by $\iota_Z :
Z \rightarrow S \tensor{S} Z$, and $\iota_Y' : Y \rightarrow Y
\tensor{S} S$ the obvious natural $S$--linear isomorphisms.

\bigskip

\noindent \emph{The functor} $-\tensor{S}\Sigma : \rcomod{S
\tensor{B} S} \rightarrow \rcomod{\rcomatrix{B}{\Sigma}}$.

\medskip Let $(Y, \rho_Y) \in \rcomod{S \tensor{B} S}$, and
consider the following right $S$--linear map
\begin{equation}\label{ros}
\xymatrix@C=60pt{Y \ar@{->}^-{\rho_Y}[r] & Y \tensfun{S}S
\tensfun{B}S \ar@{->}^-{Y\tensfun{S}\xi^{-1} \tensfun{B} S}[r] &
Y\tensfun{S} \Sigma \tensfun{A} \Sigma^* \tensfun{B} S }
\end{equation}
where $\xi$ is the $S$--bilinear map given in \eqref{xi}. Applying
$- \tensor{S}\Sigma$ to \eqref{ros}, we get
\begin{equation}\label{tensig}
\xymatrix@R=40pt@C=75pt{Y\tensfun{S} \Sigma \ar@{->}^-{\rho_Y
\tensfun{S}\Sigma}[r]\ar@{-->}_-{\rho_{Y \tensfun{S}\Sigma}}[drr]
& Y \tensfun{S} S\tensfun{B}S \tensfun{S} \Sigma \ar@{->}^-{Y
\tensfun{S} \xi^{-1} \tensfun{B} S \tensfun{S} \Sigma}[r] & Y
\tensfun{S}\Sigma \tensfun{A}\Sigma^* \tensfun{B} S \tensfun{S}
\Sigma \ar@{->}[d]|{Y \tensfun{S}\Sigma \tensfun{A} \Sigma^*
\tensfun{B} \iota^{-1}}
\\ & & Y \tensfun{S} \Sigma \tensfun{A}\Sigma^* \tensfun{B}\Sigma, }
\end{equation}
explicitly, $$\rho_{Y \tensor{S} \Sigma}(y\tensor{S}u)=
\sum_{i,(y)} y_{(0)} \tensor{S}e_i \tensor{A} e_i^* \tensor{B}
y_{(1)}u,$$ where $\rho_Y(y)= \sum_{(y)}y_{(0)} \tensor{S} 1
\tensor{B} y_{(1)}$. It is clear that $\rho_{Y \tensor{S} \Sigma}$
is a right $A$--linear map and satisfies the counitary property.
To check the coassociativity, first consider the diagram
\begin{equation}\label{diauno}
\xymatrix@R=40pt@C=90pt{Y \tensfun{S} \Sigma
\ar@{->}^-{\rho_Y\tensfun{S} \Sigma}[r]
\ar@{->}[d]|{\rho_Y\tensfun{S} \Sigma} & Y \tensfun{S} S
\tensfun{B} S \tensfun{S} \Sigma \ar@{->}[d]|{Y \tensfun{S} \Delta
\tensfun{S} \Sigma} \\ Y \tensfun{S}S \tensfun{B} S \tensfun{S}
\Sigma \ar@{->}^-{\rho_Y \tensfun{S}S \tensfun{B} S \tensfun{S}
\Sigma}[r] \ar@{->}[d]|{Y \tensfun{S} \xi^{-1}\tensfun{B} S
\tensfun{S} \Sigma} & Y \tensfun{S} S \tensfun{B} S \tensfun{S} S
\tensfun{B} S \tensfun{S} \Sigma \ar@{->}[d]|{Y \tensfun{S} S
\tensfun{B} S \tensfun{S} \xi^{-1} \tensfun{B} S \tensfun{S}
\Sigma}
\\ Y \tensfun{S} \Sigma \tensfun{A} \Sigma^* \tensfun{B} S
 \tensfun{S} \Sigma \ar@{->}^-{\rho_Y \tensfun{S} \Sigma \tensfun{A}
 \Sigma^* \tensfun{B} S \tensfun{S} \Sigma}[r] \ar@{->}[d]|{Y \tensfun{S}
  \Sigma \tensfun{A} \Sigma^* \tensfun{B} \iota^{-1}} &
  Y \tensfun{S} S \tensfun{B} S \tensfun{S} \Sigma \tensfun{A}
 \Sigma^* \tensfun{B} S
 \tensfun{S} \Sigma \ar@{->}[d]|{Y \tensfun{S} S \tensfun{B} S \tensfun{S}
 \Sigma \tensfun{A} \Sigma^* \tensfun{B} \iota^{-1}}
 \\ Y\tensfun{S} \Sigma \tensfun{A} \Sigma^* \tensfun{B}
 \Sigma \ar@{->}^-{\rho_Y\tensfun{S} \Sigma
 \tensfun{A} \Sigma^* \tensfun{B} \Sigma}[r] &
 Y \tensfun{S} S \tensfun{B} S  \tensfun{S}
 \Sigma \tensfun{A} \Sigma^* \tensfun{B}
 \Sigma }
\end{equation}
It is commutative because $\rho_Y$ is a coaction for the right $S
\tensor{B} S$--comodule $Y$. Now, look at the following diagram
\begin{equation}\label{diados}
\xymatrix@R=18pt@C=68pt{ Y \tensfun{S} S \tensfun{B} S \tensfun{S}
\Sigma \ar@{->}[dd]|{Y \tensfun{S} \Delta \tensfun{S} \Sigma}
\ar@{->}[dr]|{Y \tensfun{S} \xi^{-1} \tensfun{B}S \tensfun{S}
\Sigma} & \\ & Y \tensfun{S} \Sigma \tensfun{A} \Sigma^*
\tensor{B}S \tensor{S} \Sigma \ar@{->}[dd]|{Y \tensfun{S} \Sigma
\tensfun{A} \Sigma^* \tensfun{B} \iota^{-1}}
\\ Y \tensfun{S} S \tensfun{B} S \tensfun{S} S \tensfun{B} S
\tensfun{S} \Sigma \ar@{->}[dd]|{Y \tensfun{S} S \tensfun{B} S
\tensfun{S} \xi^{-1} \tensfun{B} S \tensfun{S} \Sigma} & \\ & Y
\tensfun{S} \Sigma \tensfun{A} \Sigma^* \tensfun{B} \Sigma
\ar@{->}[dd]|{Y \tensfun{S}\Sigma \tensfun{A} \Delta}
\\ Y \tensfun{S} S \tensfun{B} S \tensfun{S} \Sigma \tensfun{A}
 \Sigma^* \tensfun{B} S
 \tensfun{S} \Sigma \ar@{->}[dd]|{Y \tensfun{S} S \tensfun{B} S \tensfun{S}
 \Sigma \tensfun{A} \Sigma^* \tensfun{B} \iota^{-1}} &
 \\ & Y \tensfun{S}\Sigma \tensfun{A} \Sigma^* \tensfun{B} \Sigma
 \tensfun{A} \Sigma^* \tensfun{B} \Sigma^*
 \\ Y \tensfun{S} S \tensfun{B} S  \tensfun{S}
 \Sigma \tensfun{A} \Sigma^* \tensfun{B}
 \Sigma  \ar@{->}[dr]|{Y \tensfun{S} \xi^{-1}
 \tensfun{B} S \tensfun{S} \Sigma \tensfun{A} \Sigma^* \tensfun{B}
 \Sigma} &
 \\ & Y \tensfun{S} \Sigma \tensfun{A} \Sigma^*
 \tensfun{B} S \tensfun{S} \Sigma \tensfun{A} \Sigma^* \tensfun{B} \Sigma
 \ar@{->}[uu]|{Y \tensfun{S} \Sigma \tensfun{A} \Sigma^*
 \tensfun{B} \iota^{-1} \tensfun{A} \Sigma^* \tensfun{B} \Sigma}  }
\end{equation}
which is easily shown to be commutative. By concatenating diagrams
\eqref{diauno} and \eqref{diados} we see that the map
$\rho_{Y\tensor{S}\Sigma}$ endows $Y \tensor{S} \Sigma$ with a
structure of right $\rcomatrix{B}{\Sigma}$--comodule.

Now, let $f : Y \rightarrow Y'$ be a morphism in $\rcomod{S
\tensor{B} S}$, and consider the right $A$--linear map
$f\tensor{S} \Sigma : Y \tensor{S} \Sigma \rightarrow Y'
\tensor{S} \Sigma$. Then we have the following commutative diagram
\[
\xymatrix@R=20pt@C=15pt{ Y \tensfun{S} \Sigma
\ar@{->}[ddr]|{\rho_Y \tensfun{S} \Sigma} \ar@{->}[ddd]|{f
\tensfun{S} \Sigma} & & Y \tensfun{S} \Sigma \tensfun{A} \Sigma^*
\tensfun{B} S \tensfun{S} \Sigma \ar@{->}[ddr]|{Y \tensfun{S}
\Sigma \tensfun{A} \Sigma^* \tensfun{B} \iota^{-1}}
\ar@{->}[ddd]|{f\tensfun{S} \Sigma \tensfun{A} \Sigma^*
\tensfun{B} S \tensfun{S} \Sigma} & \\ & & &
\\ & Y \tensfun{S} S \tensfun{B} S \tensfun{S} \Sigma
\ar@{->}[ddd]|{f\tensfun{S} S \tensfun{B} S \tensfun{S} \Sigma}
\ar@{->}[uur]|{Y\tensfun{S} \xi^{-1} \tensfun{B} S \tensfun{S}
\Sigma} & & Y \tensfun{S} \Sigma \tensfun{A} \Sigma^* \tensfun{B}
\Sigma \ar@{->}[ddd]|{f\tensfun{S} \Sigma \tensfun{A} \Sigma^*
\tensfun{B} \Sigma}
\\ Y' \tensfun{S} \Sigma \ar@{->}[ddr]|{\rho_{Y'} \tensfun{S} \Sigma}
 & & Y'\tensfun{S} \Sigma \tensfun{A} \Sigma^* \tensfun{B}
 S \tensfun{S} \Sigma \ar@{->}[ddr]|{Y' \tensfun{S}\Sigma \tensfun{A}
 \Sigma^* \tensfun{B} \iota^{-1}}&
\\ & & &
\\ & Y' \tensfun{S} S \tensfun{B} S \tensfun{S} \Sigma
\ar@{->}[uur]|{Y' \tensfun{S} \xi^{-1} \tensfun{B} S \tensfun{S}
\Sigma} & & Y' \tensfun{S}\Sigma \tensfun{A} \Sigma^* \tensfun{B}
\Sigma, }
\]
which means that $f \tensor{S} \Sigma$ is a morphism in
$\rcomod{\rcomatrix{B}{\Sigma}}$, with the coaction
\eqref{tensig}. Therefore, we have constructed a well defined
functor $-\tensor{S} \Sigma: \rcomod{S \tensor{B} S} \rightarrow
\rcomod{\rcomatrix{B}{\Sigma}}$.

\bigskip

\noindent\emph{The functor} $-\tensor{A} \Sigma^*:
\rcomod{\rcomatrix{B}{\Sigma}} \rightarrow \rcomod{S \tensor{B}
S}$.
\medskip

 Let $(X, \rho_X) \in \rcomod{\rcomatrix{B}{\Sigma}}$, and
consider the right $S$--linear map
\begin{equation}\label{tensiges}
\xymatrix@R=40pt@C=60pt{ X\tensfun{A} \Sigma^*
\ar@{->}^-{\rho_X\tensfun{A} \Sigma^*}[r]
\ar@{-->}_-{\rho_{X\tensfun{A} \Sigma^*}}[drr] &  X
\tensfun{A}\Sigma^* \tensfun{B} \Sigma \tensfun{A} \Sigma^*
\ar@{->}^-{X \tensfun{A} \Sigma^* \tensfun{B} \xi}[r] & X
\tensfun{A} \Sigma^* \tensfun{B} S \ar@{->}[d]|{X\tensfun{A}
\iota'\tensfun{B} S}
\\ & & X\tensfun{A}\Sigma^* \tensfun{S} S \tensfun{B} S. }
\end{equation}
Direct verifications, using elements, and the coassociativity of
$\rho_X$, give a commutative diagram:
\[
\xymatrix@R=15pt@C=10pt{ X\tensfun{A} \Sigma^*
\ar@{->}[ddr]|{\rho_X \tensfun{A}\Sigma^*} \ar@{->}[dddd]|{\rho_X
\tensfun{A}\Sigma^*} & & X \tensfun{A} \Sigma^* \tensfun{B}
\ar@{->}[dddddddd]|{X \tensfun{A}\Sigma^* \tensfun{B}\mu^r} S
\\ & & \\ & X \tensfun{A} \Sigma^* \tensfun{B} \Sigma \tensfun{A} \Sigma^*
\ar@{->}[dddd]|{X \tensfun{A} \Delta \tensfun{A} \Sigma^*}
\ar@{->}[uur]|{X \tensfun{A} \Sigma^* \tensfun{B} \xi} & \\ & &
\\ X \tensfun{A} \Sigma^* \tensfun{B} \Sigma \tensfun{A} \Sigma^*
\ar@{->}[dddd]|{X \tensfun{A} \Sigma^* \tensfun{B} \xi}
\ar@{->}[ddr]|{\rho_X \tensfun{A} \Sigma^* \tensfun{B} \Sigma
\tensfun{A} \Sigma^*} & &
\\ & &
\\ & X \tensfun{A} \Sigma^* \tensfun{B} \Sigma \tensfun{A} \Sigma^*
\tensfun{B} \Sigma \tensfun{A} \Sigma^* \ar@{->}[dddd]|{X
\tensfun{A} \Sigma^* \tensfun{B} \Sigma \tensfun{A} \Sigma^*
\tensfun{B}\xi} &
\\ & &
\\ X \tensfun{A} \Sigma^* \tensfun{B} S
\ar@{->}[ddr]|{\rho_X \tensfun{A}\Sigma^* \tensfun{B} S} & & X
\tensfun{A} \Sigma^* \tensfun{B} S \tensfun{B} S
\\ & & \\ & X \tensfun{A} \Sigma^* \tensfun{B} \Sigma \tensfun{A}
\Sigma^* \tensfun{B} S \ar@{->}[uur]|{X \tensfun{A} \Sigma^*
\tensfun{B} \xi \tensfun{B} S}, & }
\]
where $\mu^r$ is the $(B-S)$--bilinear map defined by $\mu^r(s)= 1
\tensor{B} s$, for all $s \in S$. That is, the right $S$--linear
map $f:=(X \tensor{A} \Sigma^* \tensor{B} \xi ) \circ (\rho_X
\tensor{A} \Sigma^*)$ verify the cocycle condition (see
\cite[Definition 3.5(2)]{Nuss:1997}). Since $\rho_{X \tensor{A}
\Sigma^*}$ satisfies the counitary property, $f$ is actually a
descent datum on $X \tensor{A} \Sigma^*$ (see \cite{Cipolla:1976},
\cite{Nuss:1997}). Henceforth, $\rho_{X\tensor{A} \Sigma^*}= (X
\tensor{A} \iota' \tensor{B} S) \circ f$ is a right $S \tensor{B}
S$--coaction on $X \tensor{A} \Sigma^*$.

Given any right $\rcomatrix{B}{\Sigma}$--colinear map $g : X
\rightarrow X'$, we easily get a right $S \tensor{B} S$--colinear
map $g \tensor{A} \Sigma^*: X\tensor{A} \Sigma^* \rightarrow X'
\tensor{A} \Sigma^*$, with the coactions \eqref{tensiges}.
Therefore, $-\tensor{A}\Sigma^*: \rcomod{\rcomatrix{B}{\Sigma}}
\rightarrow \rcomod{S \tensor{B} S}$ is a well defined functor.

\bigskip

The precedent discussion serves to state the following
proposition.

\begin{proposition}\label{secondadj}
For every pair of comodules $\left( (Y_{S \tensor{B}S},\rho_Y);
(X_{\rcomatrix{B}{\Sigma}},\rho_X)\right)$, the following
$K$--linear map
\[
\xymatrix@R=0pt{  \Psi_{Y,X}: \hom{\rcomatrix{B}{\Sigma}}{Y
\tensfun{S}\Sigma}{X} \ar@{->}[r] & \hom{S \tensfun{B}
S}{Y}{X\tensfun{A}\Sigma^*} \\ f \ar@{|->}[r] & (f
\tensfun{A}\Sigma^*) \circ (Y\tensfun{S}\xi^{-1}) \circ \iota_Y'
\\ (X\tensfun{A}\varepsilon')
\circ (g \tensfun{S} \Sigma) & g \ar@{|->}[l]   }
\]
(where $\varepsilon'$ is the counit of the comatrix $S$--coring
$\rcomatrix{S}{\Sigma}$), is a natural isomorphism. In other
words, $-\tensor{S}\Sigma$ is left adjoint to $-\tensor{A}
\Sigma^*$.
\end{proposition}
\begin{proof}
We only prove that $\Psi_{Y,X}$ and its inverse are well defined
maps, the rest is straightforward. Clearly $\Psi_{Y,X}(f)$ is
$S$--linear, for every $f \in \hom{\rcomatrix{B}{\Sigma}}{Y
\tensor{S}\Sigma}{X}$. The colinearity of $\Psi_{Y,X}(f)$ follows
if we show that
\begin{equation}\label{dig0}
\xymatrix@R=60pt@C=60pt{ Y
\ar@{->}|-{\rho'_Y}[d] \ar@{->}^-{\Psi(f)}[r] &
X\tensfun{A}\Sigma^*  \ar@{->}|-{(X \tensfun{A} \Sigma^*
\tensfun{B} \xi)\circ (\rho_X\tensfun{A}\Sigma^*)}[d] \\
Y\tensfun{B}S \ar@{->}^-{\Psi(f)\tensfun{B}S}[r] &
X\tensfun{A}\Sigma^*\tensfun{B}S }
\end{equation}
is a commutative diagram, where $\rho'_Y=(\iota_Y^{-1}\tensor{B}S)
\circ \rho$. Put $$\boldsymbol{f}= \Psi_{Y,X}(f) \circ \rho_Y' =
(f\tensor{A} \Sigma^* \tensor{B} \xi) \circ (Y\tensor{S}\xi^{-1}
\tensor{B}S) \circ \rho_Y.$$ Using the colinearity of the map $f$,
we easily prove that the following diagram is commutative
\[
\xymatrix { Y \ar@{->}^-{\iota'}[r]
\ar@/_6pc/^{\boldsymbol{f}}[dddd] \ar@{->}[dd]|{\rho_Y} & Y
\tensfun{S} S \ar@{->}^-{Y \tensfun{S} \xi^{-1}}[r] &  Y
\tensfun{S} \Sigma \tensfun{A} \Sigma^* \ar@{->}[d]^-{\rho_Y
\tensfun{S}\Sigma \tensfun{A} \Sigma^*} \ar@/^12pc/|-{(\rho_X
\circ f)\tensfun{A}\Sigma^*}[dddd] &
\\ & & Y\tensfun{S} S \tensfun{B} S \tensfun{S} \Sigma \tensfun{A}
\Sigma^* \ar@{->}[d]|{Y \tensfun{S} \xi^{-1} \tensfun{B} S
\tensfun{S}\Sigma \tensfun{A} \Sigma^*} & \\ Y \tensfun{S} S
\tensfun{B} S \ar@{->}[d]|{Y \tensfun{S} \xi^{-1} \tensfun{B} S} &
& Y \tensfun{S} \Sigma \tensfun{A} \Sigma^* \tensfun{B} S
\tensfun{S} \Sigma \tensfun{A} \Sigma^*
\ar@{->}[d]|{Y\tensfun{S}\Sigma \tensfun{A} \Sigma^* \tensfun{B}
\iota^{-1} \tensfun{A} \Sigma^* } &
\\ Y \tensfun{S} \Sigma \tensfun{A} \Sigma^*
\tensfun{B} S \ar@{->}^-{Y \tensfun{S} \Sigma \tensfun{A} \Sigma^*
\tensfun{B} \xi^{-1}}[rr] \ar@{->}[d]|{f \tensfun{A} \Sigma^*
\tensfun{B} S} & & Y \tensfun{S} \Sigma \tensfun{A} \Sigma^*
\tensfun{B} \Sigma \tensfun{A} \Sigma^* \ar@{->}[d]|{f \tensfun{A}
\Sigma^* \tensfun{B} \Sigma \tensfun{A} \Sigma^*} & \\ X
\tensfun{A} \Sigma^* \tensfun{B} S & & X \tensfun{A} \Sigma^*
\tensfun{B} \Sigma \tensfun{A} \Sigma^* \ar@{->}_-{X \tensfun{A}
\Sigma^* \tensfun{B} \xi}[ll],  }
\]
which is exactly the diagram \eqref{dig0}. Now, let $g \in \hom{S
\tensor{B} S}{Y}{X\tensor{A}\Sigma^*}$, so the following diagram
is easily shown to be commutative
\[
\xymatrix{  Y \tensfun{S} \Sigma \ar@{->}[rr]|{g \tensfun{S}
\Sigma} \ar@{->}[ddd]|{\rho_Y \tensfun{S} \Sigma}
 & &  X \tensfun{A} \Sigma^* \tensfun{S}
\Sigma  \ar@{->}[d]|{\rho_X \tensfun{A} \Sigma^* \tensfun{S}
\Sigma}
\\ & & X\tensfun{A} \Sigma^* \tensfun{B} \Sigma \tensfun{A} \Sigma^*
\tensfun{S} \Sigma \ar@{->}[d]|{X \tensfun{A} \Sigma^* \tensfun{B}
\xi \tensfun{S} \Sigma}
\\ & & X \tensfun{A} \Sigma^* \tensfun{B}
S \tensfun{S} \Sigma \ar@{->}[d]|{X \tensfun{A} \Sigma^*
\tensfun{B}\iota^{-1}}
\\ Y \tensfun{S} S \tensfun{B} S \tensfun{S} \Sigma
\ar@{->}[dd]|{Y \tensfun{S} \xi^{-1}\tensfun{B} S \tensfun{S}
\Sigma} \ar@{->}[dr]|{\iota^{-1}_\Sigma \tensfun{B} S \tensfun{S}
\Sigma} & & X \tensfun{A} \Sigma^* \tensfun{B} \Sigma
\\ & Y \tensfun{B} S \tensfun{S} \Sigma \ar@{->}[uur]|{g \tensfun{B}
S \tensfun{S} \Sigma} &
\\ Y\tensfun{S}\Sigma
\tensfun{A} \Sigma^* \tensfun{B} S \tensfun{S} \Sigma
\ar@{->}[d]|{Y\tensfun{S}\Sigma \tensfun{A} \Sigma^*
\tensfun{B}\iota^{-1}} & & X\tensfun{A} A \tensfun{A} \Sigma^*
\tensfun{B} \Sigma \ar@{->}[uu]|{\cong}
\\ Y \tensfun{S} \Sigma \tensfun{A} \Sigma^* \tensfun{B} \Sigma
\ar@{->}[rr]|{g \tensfun{S} \Sigma \tensfun{A} \Sigma^*
\tensfun{B} \Sigma} & & X \tensfun{A} \Sigma^* \tensfun{S} \Sigma
\tensfun{A} \Sigma^* \tensfun{B}\Sigma \ar@{->}[u]|{X
\tensfun{A}\varepsilon' \tensfun{A} \Sigma^* \tensfun{B}\Sigma} }
\]
On the other hand, we have
\[
\rho_X \circ (X \tensfun{A} \varepsilon')  = (X \tensfun{A}
\Sigma^* \tensfun{B}\iota^{-1})  \circ (X \tensfun{A} \Sigma^*
\tensfun{B} \xi \tensfun{S} \Sigma)  \circ (\rho_X \tensfun{A}
\Sigma^* \tensfun{S} \Sigma),
\]
putting this in the above diagram, we get that $(X \tensor{A}
\varepsilon') \circ (g \tensor{S} \Sigma)$ is
$\rcomatrix{B}{\Sigma}$--colinear; and this finishes the proof.
\end{proof}

\begin{remark}\label{digcom}
\begin{enumerate}
\item Applying Proposition \ref{secondadj}, we get (up to natural
isomorphisms) the following commutative diagram of functors
\begin{equation}\label{dig}
\xymatrix@R=50pt@C=60pt{ & \rcomod{S \tensfun{B} S}
\ar@<0,5ex>[dr]^{\hom{S \tensfun{B} S}{S}{-}}
\ar@<0,5ex>[dl]^{-\tensfun{S}\Sigma} & \\
\rcomod{\rcomatrix{B}{\Sigma}}
\ar@<0,5ex>[rr]^{\hom{\rcomatrix{B}{\Sigma}}{\Sigma}{-}}
\ar@<0,5ex>[ur]^{-\tensfun{A}\Sigma^*} & & \rmod{B}
\ar@<0,5ex>[ll]^{-\tensfun{B}\Sigma}
\ar@<0,5ex>[ul]^{-\tensfun{B}S},}
\end{equation}
where the sideways pairs represent adjunctions. \item
Symmetrically, one can define a pair of adjoint functors relating
the categories of left comodules: $\Sigma^* \tensor{S}-: \lcomod{S
\tensor{B} S} \rightleftarrows \lcomod{\rcomatrix{B}{\Sigma}}:
\Sigma \tensor{A}-$, which turns the diagram
\begin{equation}\label{digl}
\xymatrix@R=50pt@C=60pt{ & \lcomod{S \tensfun{B} S}
\ar@<0,5ex>[dr]^{\hom{S \tensfun{B} S}{S}{-}}
\ar@<0,5ex>[dl]^{\Sigma^* \tensfun{S}-} & \\
\lcomod{\rcomatrix{B}{\Sigma}}
\ar@<0,5ex>[rr]^{\hom{\rcomatrix{B}{\Sigma}}{\Sigma^*}{-}}
\ar@<0,5ex>[ur]^{\Sigma \tensfun{A}-} & & \lmod{B}
\ar@<0,5ex>[ll]^{\Sigma^* \tensfun{B}-}
\ar@<0,5ex>[ul]^{S\tensfun{B}-},}
\end{equation}
commutative.
\end{enumerate}

\end{remark}

\section{A group isomorphism}


Let $B \subset S$ be ring extension. The set $\ibs$ of all
$B$--sub-bimodules of $S$ is a monoid with the obvious product.
For $I, J \in  \ibs$, consider the the multiplication map:
\[
\mathbf{m}: I \tensor{B} J \rightarrow IJ,\quad \mathbf{m}(x
\tensor{B} y)=xy.
\]
$\ilbs$ (resp. $\irbs$) denotes the submonoid consisting of all
$B$-sub-bimodules $I \subset S$ such that
\[
S \tensor{B} I \cong S \quad (\text{resp. } I \tensor{B} S \cong
S) \quad \text{ through } \mathbf{m}.
\]
$\invbs$ denote the group of invertible $B$-sub-bimodules of $S$.
By \cite[Proposition 1.1]{Masuoka:1989}, $\invbs \subset \ilbs
\cap \irbs$.

\bigskip
From now on fix a bimodule ${}_B\Sigma_A$ with $\Sigma_A$ finitely
generated and projective, consider its endomorphisms ring
$S=\rend{A}{\Sigma}$, and assume that ${}_B\Sigma$ is faithful, i.
e., the canonical ring extension $\lambda : B \rightarrow S$ is
injective ($B$ will be identified then with its image). Consider
the comatrix $A$--coring $\coring{C}=\rcomatrix{B}{\Sigma}$, and
denote by $\End{A}{\coring{C}}$ the monoid of the coring
endomorphisms of $\coring{C}$. We denote by $\Aut{A}{\coring{C}}$
its group of units, that is, the group of all coring automorphisms
of $\coring{C}$. The canonical Sweedler $S$--coring
$\swe{S}{B}{S}$ associated to the ring extension $B \subset S$,
will be also considered.

\begin{remark}\label{hom}
Keeping the previous notations, we made the followings remarks.
\begin{enumerate}[(1)]
\item As we have seen the $(B,A)$--bimodule $\Sigma$ is actually a
$(B,\coring{C})$--bicomodule ($B$ is considered as a trivial
$B$--coring), while $\Sigma^*$ becomes a
$(\coring{C},B)$--bicomodule. Given $g \in \End{A}{\coring{C}}$,
and a right comodule $X_{\coring{C}}$ (resp. left comodule
${}_{\coring{C}}X$), we denote by $X_g$ the associated induced
right (resp. left) $\coring{C}$--comodule. That is, $\rho_{X_g} =
(X\tensor{A}g) \, \circ \, \rho_X$ (resp. $\lambda_{X_g} =
(g\tensor{A}X) \, \circ \, \lambda_X$). If $(X,\rho_X)$ is any
right $\coring{C}$--comodule such that $X_A$ is finitely generated
and projective module, then it is well known that the right dual
module $X^*$ admits a structure of left $\coring{C}$--comodule
with coaction $$ \lambda_{X^*}(x^*) = \sum
\left((x^*\tensor{A}\coring{C}) \circ \rho_X(x_j)\right)
\tensor{A}x_j^*,\, x^* \in X^*,$$ where $\{x_j,x_j^*\}_j $ is any
right dual basis of  $X_A$.  In this way $(\Sigma_g)^*$ and
$(\Sigma^*)_g$ have the same left $\coring{C}$--coaction, that is,
they are equal as a left $\coring{C}$--comodules, then we can
remove the brackets $\Sigma_g^*=(\Sigma_g)^* = (\Sigma^*)_g$.

\item Given $g, h \in \End{A}{\coring{C}}$,  the $B$--subbimodule
$\Sigma_h \cotensor{\coring{C}} \Sigma^*_g$ of $\Sigma \tensor{A}
\Sigma^*$ is identified, via the isomorphism given in \eqref{xi},
with $\hom{\coring{C}}{\Sigma_g}{\Sigma_h}$.  Another way to
obtain this identification is given as follows. Recall, from
\cite[Example 3.4]{Gomez:2002} or \cite[Example
6]{Kaoutit/Gomez:2003a}, that $(\Sigma_g^*)_B$ is quasi-finite
$(\coring{C},B)$--bicomodule with adjunction $-\tensor{B} \Sigma_g
\dashv -\cotensor{\coring{C}} \Sigma_g^*$, so the cotensor functor
$-\cotensor{\coring{C}}\Sigma^*_g$ is naturally isomorphic to the
hom-functor $\hom{\coring{C}}{\Sigma_g}{-}$. Moreover, this
isomorphism can be chosen to be just the restriction of $-
\tensor{A} \Sigma^*_g \cong \hom{A}{\Sigma_g}{-}$. Applying this
isomorphism to $\Sigma_h$, for any $h\in \End{A}{\coring{C}}$, we
arrive to the desired identification.

\item Let $g \in\End{A}{\coring{C}}$, the following multiplication
\[
\xymatrix@R=0pt{\overline{\mathbf{m}}: \Sigma^* \tensor{B}
\hom{\coring{C}}{\Sigma_g}{\Sigma} \ar@{->}[r] & \Sigma^*_g & (u^*
\tensor{B} t \mapsto u^*t) }
\]
is a left $\coring{C}$--comodule map. Furthermore, we have a
commutative diagram
\[
\xymatrix@R=40pt@C=60pt{ \Sigma \tensfun{A} \Sigma^* \tensfun{B}
\hom{\rcomatrix{B} \Sigma}{\Sigma_g}{\Sigma} \ar@{->}^-{\Sigma
\tensfun{A} \overline{\mathbf{m}}}[r] \ar@{->}[d]|{\xi
\tensfun{A}\hom{\rcomatrix{B} \Sigma}{\Sigma_g}{\Sigma}} & \Sigma
\tensfun{A} \Sigma^*_g \ar@{->}[d]^{\xi} \\ S \tensfun{B}
\hom{\rcomatrix{B} \Sigma}{\Sigma_g}{\Sigma}
\ar@{->}^-{\mathbf{m}}[r] & S, }
\]
where $\mathbf{m}$ is the usual multiplication of $S$.
\end{enumerate}
\end{remark}

\bigskip

We define the following two maps:
\[
\xymatrix@R=0pt{ \digamma^r: \End{A}{\coring{C}} \ar@{->}[r] &
\ibs & ( g \ar@{|->}[r] & \hom{\coring{C}}{\Sigma}{\Sigma_g}), }
\]
and
\[
\xymatrix@R=0pt{ \digamma^l: \End{A}{\coring{C}} \ar@{->}[r] &
\ibs & ( g \ar@{|->}[r] & \hom{\coring{C}}{\Sigma_g}{\Sigma}). }
\]
These maps obey the following lemma. First, recall from
\cite{Sugano:1982} (cf. \cite{Caenepeel/Kadison:2001}), that $M$
is called a \emph{separable bimodule} or $B$ is said to be
$M$-\emph{separable over} $A$ provided the evaluation map $$ M
\tensor{A} {}^*M \rightarrow B,\quad m \tensor{A}\varphi \mapsto
\varphi(m)$$ is a split epimorphism of $(B,B)$--bimodules. As
shown in \cite{Sugano:1982} (cf. \cite[Theorem
3.1]{Kadison:1999}), if $M$ is a separable bimodule, then $B
\rightarrow S$ is a \emph{split extension}, i.e., there is a
$B$--linear map $\alpha: S \rightarrow B$ such that $\alpha(1_S) =
1_B$. Conversely, if ${}_BM_A$ is such that $M_A$ is finitely
generated and projective module, and $B \rightarrow S$ is a splits
extension, then ${}_BM_A$ is a separable bimodule.

\begin{lemma}\label{cos-fil}
Let $g \in \End{A}{\coring{C}}$, then
\begin{enumerate}[(i)]
\item $\digamma^r(g) \in \irbs$ if and only if ${}_BS$ preserves
the equalizer of $(\rho_{\Sigma_g} \tensor{A} \Sigma^*, \Sigma_g
\tensor{A} \lambda_{\Sigma^*})$ (cf. \cite[Section
2.4]{Gomez:2002}). In particular, if either ${}_B\Sigma$ is flat
module or ${}_B\Sigma_A$ is a separable bimodule, then
$\digamma^r(g) \in \irbs$.

\item $\digamma^l(g) \in \ilbs$ if and only if ${S}_B$ preserves
the equalizer of $(\rho_{\Sigma} \tensor{A} \Sigma_g^*, \Sigma
\tensor{A} \lambda_{\Sigma_g^*})$. In particular, if either
$\Sigma_B^*$ is flat module or ${}_B\Sigma_A$ is a separable
bimodule, then $\digamma^l(g) \in \ilbs$.

\item If $g \in \Aut{A}{\coring{C}}$, then $\digamma^l(g) =
\digamma^r(g^{-1})$.
\end{enumerate}
\end{lemma}
\begin{proof}
$(i)$ and $(ii)$ We only prove $(i)$ because $(ii)$ is symmetric.
Following the identifications made in Remark \ref{hom}, we have
$\digamma^r(g) \cong \Sigma_g\cotensor{\coring{C}}\Sigma^*$.
Taking this isomorphism into account, the first statement in $(i)$
is reduced to the problem of compatibility between tensor and
cotensor. Effectively, by \cite[Lemma 2.2]{Gomez:2002}, ${}_BS
\cong \Sigma \tensor{A} \Sigma^*$ preserves the equalizer of
$(\rho_{\Sigma_g} \tensor{A} \Sigma^*, \Sigma_g \tensor{A}
\lambda_{\Sigma^*})$ if and only if
\[
(\Sigma_g \cotensor{\coring{C}} \Sigma^*) \tensor{B} \Sigma
\tensor{A} \Sigma^*  \cong \Sigma_g \cotensor{\coring{C}}
(\Sigma^* \tensor{B} \Sigma \tensor{A} \Sigma^*) = (\Sigma_g
\cotensor{\coring{C}} \coring{C}) \tensor{A} \Sigma^* \cong
\Sigma\tensor{A} \Sigma^* \cong S,
\]
if and only if $(\Sigma_g \cotensor{\coring{C}} \Sigma^*) \in
\irbs$, since by Remark \ref{hom}(3) this composition coincides
with the multiplication of the monoid $\ibs$. If ${}_B\Sigma$ is a
flat module, then clearly ${}_BS$ is also flat. Hence, it
preserves the stated equalizer. Now, if we assume that
${}_B\Sigma_A$ is a separable bimodule, then \cite[Theorem
3.5]{Brzezinski/Gomez:2003} implies that
$\coring{C}=\rcomatrix{B}{\Sigma}$ is a coseparable $A$--coring
(cf. \cite{Guzman:1989}, \cite{Gomez/Louly:2003} for definition).
Therefore, equalizers split by \cite[Proposition
1.2]{Guzman:1989}, and so they are preserved by any module.\\
$(iii)$ A straightforward computation shows that
$\hom{\coring{C}}{\Sigma_g}{\Sigma} =
\hom{\coring{C}}{\Sigma}{\Sigma_{g^{-1}}}$.
\end{proof}

\begin{theorem}\label{resultado-0}
Let ${}_B \Sigma_A$ be a bimodule such that ${}_B\Sigma$ is
faithful and $\Sigma_A$ is finitely generated and projective.
Consider $\coring{C}=\rcomatrix{B}{\Sigma}$ its associated
comatrix $A$--coring. If either
\begin{enumerate}[(a)]
\item $\Sigma^*_B$ is a faithfully flat module, or \item
${}_B\Sigma_A$ is a separable bimodule.
\end{enumerate}
Then $\digamma^l: \End{A}{\coring{C}} \rightarrow \ilbs$ is a
monoid isomorphism with inverse
\begin{equation}\label{Gama}
\xymatrix@R=0pt{\Gamma^l:\ilbs \ar@{->}[r] & \End{A}{\coring{C}}
\\ I \ar@{|->}[r] & [u^* \tensfun{B} u \mapsto \sum_k u^*
s_k \tensfun{B}  x_ku], }
\end{equation}
where $\mathbf{m}^{-1}(1)=\sum_k s_k \tensor{B} x_k \in
S\tensor{B}I$.
\end{theorem}
\begin{proof}
Under the hypothesis $(a)$, we have, by the left version of the
generalized Descent Theorem for modules \cite[Theorem
2]{Kaoutit/Gomez:2003a}, that $\Sigma^* \tensor{B}-: \lmod{B}
\rightarrow \lcomod{\rcomatrix{B}{\Sigma}}$ is an equivalence of
categories with inverse
$\hom{\rcomatrix{B}{\Sigma}}{\Sigma^*}{-}$. Applying the diagram
\eqref{digl} of the Remark \ref{digcom}, we obtain that $S
\tensor{B}-: \lmod{B} \rightarrow \lcomod{S \tensor{B}S}$ is a
separable functor (cf. \cite{Nastasescu/Bergh/Oystaeyen:1989} for
definition). Now, assume $(b)$, then the ring extension $B
\rightarrow S$ splits as a $B$--bimodule map. By \cite[Proposition
1.3]{Nastasescu/Bergh/Oystaeyen:1989}, the functors $S\tensor{B}-:
\lmod{B} \rightarrow \lmod{S}$ is separable, and by \cite[Lemma
1.1(3)]{Nastasescu/Bergh/Oystaeyen:1989}, the functor
$S\tensor{B}-: \lmod{B} \rightarrow \lcomod{S\tensor{B}S}$ is
separable. In conclusion, under the hypothesis $(a)$ or $(b)$, the
functor $S\tensor{B}-: \lmod{B} \rightarrow \lcomod{S\tensor{B}S}$
reflects isomorphisms. Therefore, any inclusion $I \subseteq J$ in
$\ilbs$, implies equality $I=J$. This fact will be used implicitly
in the remainder of the proof.

The map $\Gamma^l$ is easily shown to be well defined, while Lemma
\ref{cos-fil} implies that $\digamma^l$ is also well defined. Let
us first show that $\digamma^l$ is a monoid map. The image of the
unit is mapped to $B$, $\digamma^l(1_{\End{A}{\coring{C}}})=
\mathrm{End}_{\coring{C}}(\Sigma)= B$, since by \cite[Proposition
2]{Kaoutit/Gomez:2003a} the inclusion $B \subseteq
\rend{\rcomatrix{B}{\Sigma}}{\Sigma}$ is always true. Let $g,h \in
\End{A}{\coring{C}}$, and $t \in \digamma^l(g)$, $s \in
\digamma^l(h)$, that is \begin{eqnarray*}
  \sum_i e_i \tensor{A} e_i^* \tensor{B} tu &=& \sum_i te_i\tensor{A}g(e_i^* \tensor{B}u ) \\
 \sum_i e_i \tensor{A} e_i^* \tensor{B} su  &=& \sum_i se_i\tensor{A}h(e_i^* \tensor{B}u )
\end{eqnarray*}
for every element $u \in \Sigma$. So, for every $u \in \Sigma$, we
have
\begin{eqnarray*}
  \rho_{\Sigma}(tsu) &=& \sum_i e_i\tensor{A}e_i^* \tensor{B}tsu  \\
   &=& \sum_i te_i\tensor{A}g(e_i^* \tensor{B}su ) \\
   &=& (t\tensor{A}\coring{C}) \, \circ \, (\Sigma \tensor{A} g)
   \left( \sum_i e_i\tensor{A}e_i^* \tensor{B}su \right) \\
   &=& (t\tensor{A}\coring{C}) \, \circ \, (\Sigma \tensor{A} g)
   \left( \sum_i se_i\tensor{A}h(e_i^* \tensor{B}u) \right)  \\
   &=& \sum_i tse_i \tensor{A} gh(e_i^*\tensor{B}u) \\
   &=& (ts\tensor{A}\coring{C}) \, \circ \, \rho_{\Sigma_{gh}}(u)
\end{eqnarray*}
which means that $ts \in
\hom{\coring{C}}{\Sigma_{gh}}{\Sigma}=\digamma^l(gh)$, and so
$\digamma^l(g) \digamma^l(h) = \digamma^l(gh)$. Now, let $I \in
\ilbs$ with $\mathbf{m}^{-1}(1)= \sum_k s_k \tensor{B} t_k \in S
\tensor{B} I$. If $s$ is any element in $I$, then $1\tensor{B}s =
\sum_kss_k \tensor{B}t_k \in S\tensor{B}I$. Henceforth, \begin{eqnarray*}
  (s\tensor{A}\coring{C}) \, \circ \, \rho_{\Sigma_{\Gamma^l(I)}}(u) &=&
  (s\tensor{A}\coring{C})
  \left( \sum_i e_i\tensor{A}\Gamma^l(I)(e_i^*\tensor{B}u) \right) \\
   &=& \sum_{i,k} se_i \tensor{A}e_i^* s_k \tensor{B}t_ku \\
   &=& \sum_{i,k} e_i \tensor{A}e_i^* s s_k \tensor{B}t_ku \\
   &=& \sum_{i} e_i \tensor{A}e_i^*  \tensor{B}su \, = \,
   \rho_{\Sigma}(su)
\end{eqnarray*}
for every $u \in \Sigma$, that is $s: \Sigma_{\Gamma^l(I)}
\rightarrow \Sigma \in I$ is a $\coring{C}$--colinear map.
Therefore, $I = \digamma^l(\Gamma^l(I))$, for every $I \in \ilbs$.
Conversely, let $g \in \End{A}{\rcomatrix{B}{\Sigma}}$, and put
$I=\digamma^l(g)= \hom{\coring{C}}{\Sigma_g}{\Sigma}$ with
$\mathbf{m}^{-1}(1)=\sum_k s_k \tensor{B} x_k \in S\tensor{B}I$.
For every $t \in I$, we have
\begin{equation}\label{uphi}
\begin{array}{ll}
\sum_i g(u^* t \tensor{B} e_i) \tensor{A} e_i^* = \sum_i u^*
\tensor{B} e_i \tensor{A} e_i^* t, & \quad \forall u^* \in
\Sigma^*
\\ & \\ \sum_i e_i \tensor{A} e_i^* \tensor{B}  tu = \sum_i te_i
\tensor{A} g(e_i^* \tensor{B} u), & \quad \forall u \in \Sigma.
\end{array}
\end{equation}
Computing, using equations \eqref{uphi}
\begin{eqnarray*}
  (\Gamma^l(I)\tensor{A} \lambda_{\Sigma^*}(u^*) &=& \sum_{i,k} u^*s_k \tensor{B} t_ke_i \tensor{A} e_i^* \\
   &=& \sum_{k} u^*s_k \tensor{B} \left( \sum_i t_ke_i \tensor{A} e_i^* \right) \\
   &=& \sum_{k} u^*s_k \tensor{B} \left( \sum_i e_i \tensor{A} e_i^* t_k\right) \\
   &=& \sum_k\left( \sum_i u^*s_k \tensor{B} e_i \tensor{A}e_i^*t_k \right) \\
   &=& \sum_k\left( \sum_i g(u^*s_kt_k \tensor{B} e_i) \tensor{A}e_i^* \right) \\
   &=& \sum_i g(u^*\tensor{B}e_i)\tensor{A}e_i^* \\
   &=& (g \tensor{A} \Sigma^*) \, \circ \, \lambda_{\Sigma^*}(u^*)
\end{eqnarray*} for every $u^* \in \Sigma^*$, that is $(\Gamma(I) \tensor{A} \Sigma^*) \circ
\lambda_{\Sigma^*} = (g \tensor{A} \Sigma^*) \circ
\lambda_{\Sigma^*}$. Whence,
\begin{equation}\label{delt-g}
(\Gamma(I) \tensor{A} \Sigma^* \tensor{B} \Sigma) \circ \Delta =
(g \tensor{A} \Sigma^* \tensor{B} \Sigma) \circ \Delta,
\end{equation}
because $\Delta_{\coring{C}} =  \lambda_{\Sigma^*} \tensor{B}
\Sigma$. On the other hand
\begin{eqnarray*}
  \Delta \, \circ \, \Gamma^l(I)(u^* \tensor{B} u) &=& \sum_k u^* s_k
\tensor{B}\left(\sum_i e_i \tensor{A} e_i^* \tensor{B} t_ku\right) \\
   &=& \sum_k u^* s_k \tensor{B} \left(\sum_i t_ke_i \tensor{A} g(e_i^*
\tensor{B} u)\right), \, \text{ by } \eqref{uphi} \\
   &=& \Sigma^* \tensor{B} \Sigma \tensor{A} g\left(\sum_{i,k} u^* s_k
\tensor{B}t_ke_i \tensor{A} e_i^* \tensor{B} u\right) \\
   &=& (\Sigma^* \tensor{B} \Sigma \tensor{A} g) \, \circ \, (\Gamma^l(I)
\tensor{A} \Sigma^* \tensor{B} \Sigma) \, \circ \, \Delta (u^* \tensor{B}u) \\
   &=& (\Sigma^* \tensor{B} \Sigma \tensor{A} g) \,  \circ \, (g
\tensor{A} \Sigma^* \tensor{B} \Sigma) \, \circ \, \Delta(u^*
\tensor{B}u),\, \text{ by } \eqref{delt-g} \\
   &=& (g \tensor{A} g) \, \circ \, \Delta (u^* \tensor{B} u) \\
   &=& \Delta \, \circ \, g(u^* \tensor{B}u), \, g \in
   \End{A}{\coring{C}},
\end{eqnarray*}
for every $u^*  \in \Sigma^*,\, u \in \Sigma$. Therefore, $\Delta
\circ \Gamma^l(I) = \Delta \circ g$, thus
$\Gamma^l(I)=\Gamma^l(\digamma^l(g))=g$, for every $g \in
   \End{A}{\coring{C}}$ since $\Delta$
is injective.
\end{proof}

Symmetrically we have the anti-homomorphism of monoids
\begin{equation}\label{Gama'}
\xymatrix@R=0pt{ \Gamma^r: \irbs \ar@{->}[r] &
\End{A}{\rcomatrix{B}{\Sigma}} \\ I \ar@{|->}[r] & [u^* \tensor{B}
u \mapsto \sum_k u^* t_k \tensfun{B} s_k u], }
\end{equation}
where $\mathbf{m}^{-1}(1) = \sum_k t_k \tensor{B} s_k \in I
\tensor{B} S$. Let $B^o \subset S^o$ denote the opposite ring
extension of $B \subset S$, and identify $S^o$ with
$\rend{A^o}{(\Sigma^*)^o}$, where the notation $X^o$, for any left
$A$--module $X$, means the opposite right $A^o$--module. Put
${}_{B^o} W_{A^o}= ({}_A \Sigma^*_B)^o$ the opposite bimodule, and
consider its right dual $W^*$, with respect to $A^o$, i.e. $W^*
=\mathrm{Hom}(W_{A^o},A^o_{A^o})$. Obviously $W_{A^o}$ is finitely
generated and projective module, and we can consider its
associated comatrix $A^o$--coring $\rcomatrix{B^o}{W}$. By the
Remark \ref{lcomcor}, there is an $A$--coring isomorphism
\[
(W^* \tensor{B^o} W)^o \cong \rcomatrix{B}{\Sigma}, \qquad \left(
(w^* \tensor{B^o} w)^o \mapsto \sum_i w \tensor{B} e_i
w^*((e_i^*)^o)^o \right),
\]
where $(W^* \tensor{B^o} W)^o$ is the opposite $A$--coring of the
$A^o$--coring $W^* \tensor{B^o} W$. Therefore, we have an
isomorphism of monoids $\End{A^o}{\rcomatrix{B^o}{W}} \cong
\End{A}{\rcomatrix{B}{\Sigma}}$. Finally, using this last
isomorphism together with the equality $\irbs =
\mathbf{I}^l_{B^o}(S^o)$, we can identify the $\Gamma^r$-map of
equation \eqref{Gama'} with the $\Gamma^l$-map \eqref{Gama}
associated to the new data: $A^o$, $B^o \subset S^o$, and
${}_{B^o}W_{A^o}$. Henceforth, Theorem \ref{resultado-0} yields

\begin{theorem}\label{resultado-0'}
Let ${}_B \Sigma_A$ be a bimodule such that ${}_B\Sigma$ is
faithful and $\Sigma_A$ is finitely generated and projective.
Consider $\coring{C}=\rcomatrix{B}{\Sigma}$ its associated
comatrix $A$--coring. If either
\begin{enumerate}[(a)]
\item ${}_B\Sigma$ is faithfully flat module, or \item
${}_B\Sigma_A$ is a separable bimodule.
\end{enumerate}
Then $\digamma^r: \End{A}{\rcomatrix{B}{\Sigma}} \rightarrow
\irbs$ is an anti-isomorphism of monoids with inverse map
$$\xymatrix@R=0pt{
\Gamma^r: \irbs \ar@{->}[r] & \End{A}{\rcomatrix{B}{\Sigma}} \\ I
\ar@{|->}[r] & [u^* \tensor{B} u \mapsto \sum_k u^* t_k
\tensfun{B} s_k u], }$$ where $\mathbf{m}^{-1}(1) = \sum_k t_k
\tensor{B} s_k \in I \tensor{B} S$.
\end{theorem}

The isomorphism $\Gamma^l$ given in \eqref{Gama} gives, by
restriction, an isomorphism of groups $\Gamma : \invbs \rightarrow
\Aut{A}{\rcomatrix{B}{\Sigma}}$. Analogously, the anti-isomorphism
$\Gamma^r$ defined in \eqref{Gama'}, gives, by restriction, an
anti-isomorphism of groups $\Gamma' : \invbs \rightarrow
\Aut{A}{\rcomatrix{B}{\Sigma}}$. Moreover, when both $\Gamma^r$
and $\Gamma^l$ are bijective, Lemma \ref{cos-fil}.(iii) says that
$\Gamma = (-)^{-1} \circ \Gamma'$, where $(-)^{-1}$ denotes the
antipode map in the group of automorphisms. We can thus say that,
either in the hypotheses of Theorem \ref{resultado-0} or in the
hypotheses of Theorem \ref{resultado-0'}, we have an isomorphism
of groups $\Gamma : \invbs \rightarrow
\Aut{A}{\rcomatrix{B}{\Sigma}}$ defined either as $\Gamma^l$ or as
$(-)^{-1} \circ \Gamma^r$, respectively. We can then state our
main theorem as follows.

\begin{theorem}\label{resultado-1}
Let ${}_B \Sigma_A$ be a bimodule such that ${}_B\Sigma$ is
faithful and $\Sigma_A$ is finitely generated and projective.
Consider $\coring{C}=\rcomatrix{B}{\Sigma}$ its associated
comatrix $A$--coring. If either
\begin{enumerate}[(a)]
\item ${}_B\Sigma$ or $\Sigma^*_B$ is a faithfully flat module, or
\item ${}_B\Sigma_A$ is a separable bimodule.
\end{enumerate}
Then there is an isomorphism of groups $\Gamma: \invbs \rightarrow
\Aut{A}{\rcomatrix{B}{\Sigma}}$.
\end{theorem}

To finish, we want to compare Masuoka's maps \cite[Theorem
2.2(2.3)]{Masuoka:1989} with our $\digamma$--maps, using the
adjunction of the section \ref{Sect1}.

\begin{proposition}\label{(b)1}
Let ${}_B \Sigma_A$ be a bimodule such that ${}_B\Sigma$ is
faithful and $\Sigma_A$ is finitely generated and projective. Let
$S=\rend{A}{\Sigma}$ its ring of right linear endomorphisms. Then
\begin{enumerate}[(1)]
\item the map
\[
\xymatrix@R=0pt{\widehat{(-)}: \End{A}{\rcomatrix{B}{\Sigma}}
\ar@{->}[r] & \End{S}{S \tensor{B} S} \\ g \ar@{|->}[r] &
\widehat{g}= (\xi \tensor{B} \xi) \circ (\Sigma \tensor{A} g
\tensor{A}\Sigma^*) \circ (\xi^{-1}  \tensor{B} \xi^{-1}) }
\]
is an injective homomorphism of monoids which turns the following
diagram commutative
\[
\xymatrix{ \ilbs \ar[d]_{\overline{\Gamma}^l} \ar[r]^-{\Gamma^l} &
\End{A}{\rcomatrix{B}{\Sigma}} \ar[dl]^-{\widehat{(-)}}
\\
\End{S}{S \tensfun{B} S} &  }
\]
where $\overline{\Gamma}^l$ is the Gamma map associated to the
bimodule ${}_BS_S$ and the comatrix $S$--coring $S\tensor{B}S$
(see \cite[(2.1)]{Masuoka:1989}); \item for every $g \in
\End{A}{\rcomatrix{B}{\Sigma}}$, we have
\[
\hom{\rcomatrix{B}{\Sigma}}{\Sigma_g}{\Sigma} \,=\,
\hom{S\tensor{B}S}{S_{\widehat{g}}}{S} \, = \, \{s \in S |\,\,
\widehat{g}(s\tensor{B}1)=1 \tensor{B}s \}
\]
\end{enumerate}
\end{proposition}
\begin{proof}
(1) We only show that $\widehat{(-)}$ is a well defined map, the
compatibilities with the multiplication and unit are an easy
computations. So let $g \in \End{A}{\rcomatrix{B}{\Sigma}}$, by
definition $\widehat{g}$ is an $S$--bilinear map, and preserves
the counit. Denote by $\Delta'$ the comultiplication of $S
\tensor{B} S$, i.e. $\Delta': S\tensor{B}S \rightarrow
S\tensor{B}S\tensor{B}S$ sending $s\tensor{B}s' \mapsto s
\tensor{B}1\tensor{B}s'$, $s, s'\in S$. Then $\widehat{g}$ is
coassociative if and only if
\begin{equation}\label{coass}
\Delta' \circ \widehat{g} = (\widehat{g} \tensfun{B} S) \circ (S
\tensfun{B} \widehat{g}) \circ \Delta'.
\end{equation}
Now, a direct computations give the following equations
\begin{eqnarray*}
  ( \widehat{g} \tensor{B} S) \circ (S \tensor{B} \widehat{g}) &=&
  (\xi \tensor{B} \xi \tensor{B} \xi) \circ ( \Sigma \tensor{A} g
\tensor{A} g \tensor{A} \Sigma) \circ (\xi^{-1} \tensor{B}
\xi^{-1} \tensor{B} \xi^{-1}),
\end{eqnarray*}
\begin{eqnarray*}
(\Sigma \tensor{A} \Delta \tensor{A} \Sigma^* ) \circ (\xi^{-1}
\tensor{B} \xi^{-1}) &=& (\xi^{-1} \tensor{B} \xi^{-1} \tensor{B}
\xi^{-1}) \circ \Delta, \\ & & \\
 \Delta'  \circ (\xi \tensor{B} \xi)  &=& (\xi \tensor{B} \xi \tensor{B} \xi) \circ (\Sigma \tensor{A}
\Delta \tensor{A} \Sigma^*),
\end{eqnarray*}
which in conjunction with the coassociativity of $g$ imply the
equality of equation \eqref{coass}.
\\ (2) The second stated equality is a direct consequence of
the identification of the $B$--bimodule
$\hom{S\tensor{B}S}{S_{\widehat{g}}}{S}$ with a $B$--sub-bimodule
of $S$. Now, observe that the canonicals right $A$--linear and
right $S$--linear isomorphisms $S_{\widehat{g}} \tensor{S} \Sigma
\cong \Sigma_g$ and $S \cong \Sigma \tensor{A} \Sigma^*$ are,
respectively, right $\rcomatrix{B}{\Sigma}$--colinear map and
right $S \tensor{B} S$--colinear map, with respect to the
coactions defined in equations \eqref{tensig} and
\eqref{tensiges}. Whence,
\[
\hom{\rcomatrix{B}{\Sigma}}{\Sigma_g}{\Sigma} \cong
\hom{\rcomatrix{B}{\Sigma}}{S_{\widehat{g}}\tensor{S} \Sigma
}{\Sigma} \cong \hom{S \tensor{B} S}{S_{\widehat{g}}}{\Sigma
\tensor{A} \Sigma^*} \cong \hom{S \tensor{B}
S}{S_{\widehat{g}}}{S},
\]
where the second isomorphism is given by the Proposition
\ref{secondadj}. The desired first equality is now obtained using
the inclusion $\hom{\rcomatrix{B}{\Sigma}}{\Sigma_g}{\Sigma}
\subseteq \hom{S \tensor{B} S}{S_{\widehat{g}}}{S} \subset S$
which we show as follows. An element $s \in S$ belongs to
$\hom{\rcomatrix{B}{\Sigma}}{\Sigma_g}{\Sigma}$ if and only if
\begin{eqnarray*}
  \sum_i e_i \tensor{A} e_i^* \tensor{B}  su &=& \sum_i se_i
\tensor{A} g(e_i^* \tensor{B} u), \quad \forall u \in \Sigma.
\end{eqnarray*} This implies  \begin{eqnarray*}
  \sum_{i,j} e_i \tensor{A} e_i^* \tensor{B}  se_j \tensor{A} e_j^* &=& \sum_{i,j} se_i
\tensor{A} g(e_i^* \tensor{B} e_j)\tensor{A}e_j^*
\end{eqnarray*} Using the isomorphism $\xi$ of equation
\eqref{xi} and the definition of the map $\widehat{(-)}$, we
obtain $s \in \hom{\rcomatrix{B}{\Sigma}}{\Sigma_g}{\Sigma}$
implies $1\tensor{B}s = \widehat{g}(s\tensor{B}1)$.
\end{proof}

\providecommand{\bysame}{\leavevmode\hbox
to3em{\hrulefill}\thinspace}
\providecommand{\MR}{\relax\ifhmode\unskip\space\fi MR }
\providecommand{\MRhref}[2]{%
  \href{http://www.ams.org/mathscinet-getitem?mr=#1}{#2}
} \providecommand{\href}[2]{#2}

\end{document}